\documentclass{gtpart}

\usepackage{graphicx,psfrag}

\title{Saddle tangencies and the distance of Heegaard splittings}
\author{Tao Li}
\givenname{Tao}
\surname{Li}

\address{Department of Mathematics \\
Boston College \\
Chestnut Hill, MA 02467\\USA}
\email{taoli@bc.edu}
\urladdr{http://www2.bc.edu/~taoli/}

\keyword{Heegaard splitting}
\keyword{incompressible surface}
\keyword{curve complex}
\keyword{sample layout}
\subject{primary}{msc2000}{57N10}
\subject{secondary}{msc2000}{57M50}

\input xypic
\let\xysavmatrix\xymatrix
\def\xymatrix{\disablesubscriptcorrection\xysavmatrix}
\AtBeginDocument{\let\bar\wbar\let\hat\what}
\def\figdir{.}

\makeop{genus}
\newcommand{\Int}{\mathrm{int}}

\theoremstyle{plain}
\newtheorem{theorem}{Theorem}[section]
\newtheorem{lemma}{Lemma}[section]
\newtheorem{corollary}{Corollary}[section]

\theoremstyle{definition}
\newtheorem{definition}{Definition}[section]

\theoremstyle{remark}
\newtheorem*{remark}{Remark}

\newtheorem{claim}{Claim}
\newtheorem*{claima}{Claim A}
\newtheorem*{claimb}{Claim B}
\newtheorem*{labelling}{Labelling}

\newtheorem*{notation}{Notation}

\numberwithin{figure}{section}

\begin{document}
\begin{psfrags}

\makeatletter
\let\c@corollary\c@lemma\c@theorem
0\makeatother

\begin{abstract}
We give another proof of a theorem of Scharlemann and Tomova and of a
theorem of Hartshorn. The two theorems together say the following.
Let $M$ be a compact orientable irreducible 3--manifold and $P$ a
Heegaard surface of $M$.  Suppose $Q$ is either  an incompressible
surface or a strongly irreducible Heegaard surface in $M$. Then either
the Hempel distance $d(P)\le 2\genus(Q)$ or $P$ is isotopic to $Q$.  This theorem can be naturally extended to bicompressible
but weakly incompressible surfaces.
\end{abstract}

\maketitle

\section{Introduction}\label{Sintro}
Let $P$ be a closed orientable surface of genus at least 2.  The curve
complex of $P$, introduced by Harvey \cite{H}, is the complex whose
vertices are the isotopy classes of essential simple closed curves in
$P$, and $k+1$ vertices determine a $k$--simplex if they are
represented by pairwise disjoint curves.  We denote the curve complex
of $P$ by $\mathcal{C}(P)$.  For any two vertices in $\mathcal{C}(P)$,
the distance $d(x,y)$ is the minimal number of 1--simplices in a
simplicial path jointing $x$ to $y$.  To simplify notation, unless
necessary, we do not distinguish a vertex in $\mathcal{C}(P)$ from a
simple closed curve in $P$ representing this vertex.

Let $M$ be a compact orientable irreducible  3--manifold and $P$ an
embedded connected separating surface in $M$ with $\genus(P)\ge 2$.
Let $U$ and $V$ be the closure of the two components of $M-P$. We may
view $\partial U=\partial V=P$. As in Scharlemann--Tomova
\cite{ST}, we say $P$ is
\emph{bicompressible\/} if $P$ is compressible in both $U$ and $V$.  Let
$\mathcal{U}$ and $\mathcal{V}$ be the set of vertices in
$\mathcal{C}(P)$ represented by curves bounding compressing disks in
$U$ and $V$ respectively.  The distance $d(P)$ is defined to be the
distance between $\mathcal{U}$ and $\mathcal{V}$ in the curve complex
$\mathcal{C}(P)$.  If $P$ is a Heegaard surface, then $d(P)$ is the
distance defined by Hempel \cite{He}.  We say $P$ is \emph{strongly
irreducible\/} or following the definition in \cite{ST}, say $P$ is
\emph{weakly incompressible\/} if $d(P)\ge 2$, ie every compressing
disk in $U$ intersects every compressing disk in $V$.

Let $Q$ be another closed orientable surface embedded in $M$.  Let
$g(Q)$ be the genus of $Q$.  A theorem of Hartshorn \cite{Ha} says
that if $Q$ is incompressible and $P$ is a strongly irreducible
Heegaard surface, then $d(P)\le 2g(Q)$.  In \cite{ST}, Scharlemann and
Tomova showed that if both $P$ and $Q$ are connected, separating, bicompressible
and strongly irreducible, then either $d(P)\le 2g(Q)$ or $P$ and $Q$
are well-separated or $P$ and $Q$ are isotopic.  In particular, if both $P$ and $Q$
are strongly irreducible Heegaard surfaces, either $P$ and $Q$ are
isotopic or $d(P)\le 2g(Q)$.

Combining Hartshorn's theorem and the theorem of Scharlemann and
Tomova, we have the following Theorem.

\begin{theorem}\label{Tmain}
Suppose $M$ is a compact orientable irreducible  3--manifold and $P$
is a separating bicompressible and strongly irreducible (or weakly
incompressible) surface in $M$.  Let $Q$ be an embedded closed
orientable surface in $M$ and suppose $Q$ is either incompressible or
separating, bicompressible but strongly irreducible.  Then either
\begin{enumerate}
\item $d(P)\le 2g(Q)$, or
\item after isotopy, $P_t\cap Q=\emptyset$ for all $t$, where $P_t$ ($t\in [0,1]$) is a level surface in a sweep-out for $P$, see \fullref{Skey} for definition, or
\item $P$ and $Q$ are isotopic.
\end{enumerate}
\end{theorem}

\begin{remark} The statement of \fullref{Tmain} is basically the
same as  the main theorem of \cite{ST}. If $Q$ is separating, bicompressible but strongly irreducible and $P_t\cap Q=\emptyset$ for all $t\in [0,1]$, then it is easy to see that $P$ and $Q$ are well-separated. Note that part (3) of the theorem never happens if $Q$ is incompressible.
\end{remark}

In this paper, we give another proof of \fullref{Tmain}.  Some
arguments were originally used in a different proof of the main
theorem by the author \cite{L1}.  The motivation for this paper is a conjecture
in \cite{L1} which generalizes both the main theorem of \cite{L1} and
the theorem of Scharlemann and Tomova.  We hope this proof and the
techniques in \cite{L1,L2} can lead to a solution of this conjecture.
Some arguments in the proof are similar to those in
\cite{BS, ST}.  

I would like to thank Marty Scharlemann for pointing out a mistake in an earlier version of the paper. The research was partially supported by NSF grant DMS-0406038.

\section{Saddle tangencies}\label{Skey}

\begin{notation}
Throughout this paper, we denote the interior of $X$ by $\Int(X)$ for
any space $X$.
\end{notation}

Let $P$ be a bicompressible surface and let $U$ and $V$ be the closure
of the two components of $M-P$ as above.  Let $P^U$ and $P^V$ be the
possibly disconnected surfaces obtained by maximal compressing $P$ in
$U$ and $V$ respectively.  Since $M$ is irreducible, after capping off
2-sphere components by 3-balls, we may assume $P^U$ and $P^V$ do not
contain 2-sphere components.  Moreover, we may also assume $P^U\subset
\Int(U)$ and $P^V\subset \Int(V)$.  Since $P$ is strongly irreducible,
as in Casson--Gordon 
\cite{CG}, $P^U$ and $P^V$ are incompressible in $M$.
Furthermore, $P^U\cup P^V$ bounds a submanifold $M_P$ of $M$ and $P$
is a strongly irreducible Heegaard surface of $M_P$.  Note that if $U$
is a handlebody, then $P^U=\emptyset$.  If $P$ is a Heegaard surface
of $M$, then we may view $M_P=M$.

The surface $P$ cuts $M_P$ into a pair of compression bodies $U\cap
M_P$ and $V\cap M_P$.  There are a pair of properly embedded graphs
$G^U\subset U\cap M_P$ and $G^V\subset V\cap M_P$ which are the spines
of the two compression bodies.  The endpoints of the graphs $G^U$ and
$G^V$ lie in $P^U$ and $P^V$ respectively.  Let $\Sigma_U=P^U\cup G^U$
and $\Sigma_V=P^V\cup G^V$, then $M_P-(\Sigma_U\cup\Sigma_V)$ is
homeomorphic to $P\times (0,1)$.

We consider a sweepout $H \co P\times (I,\partial I)\to
(M_P,\Sigma_U\cup\Sigma_V)$, see \cite{RS}, where $I=[0,1]$ and
$H|_{P\times(0,1)}$ is an embedding.  We denote $H(P\times\{x\})$ by
$P_x$ for any $x\in I$. We may assume $P_0=\Sigma_U$, $P_1=\Sigma_V$
and each $P_x$ ($i\ne 0,1$) is isotopic to $P$.  To simplify notation,
we will not distinguish $H(P\times (0,1))$ from $P\times (0,1)$.

Let $\pi \co P\times I\to P$ be the projection.  To simplify notation, we
do not distinguish between an essential simple closed curve $\gamma$
in $P_x$ and the vertex represented by $\pi(\gamma)$ in the curve
complex $\mathcal{C}(P)$.

\begin{definition}
Let $Q$ be a properly embedded compact surface in $M$.  We say $Q$ is
in \emph{regular position with respect to $P\times I$\/} if
\begin{enumerate}
\item $Q\cap G^U$ and $Q\cap G^V$ consist of finitely many points and
$Q$ is transverse to $P^U\cup P^V$ and
\item $Q$ is transverse to each $P_x$, $x\in (0,1)$, except for
finitely many critical levels $t_1,\dots, t_n\in (0,1)$ and
\item at each critical level $t_i$, $Q$ is transverse to $P_{t_i}$
except for a single saddle or center tangency.
\end{enumerate}
If $x\in (0,1)$ is not one of the $t_i$'s, then we say $x$ or $P_x$ is
a regular level.  Clearly every embedded surface $Q$ can be isotoped
into a regular position.
\end{definition}

\begin{definition}
We say $Q$ is \emph{irreducible with respect to $P\times I$\/} if
\begin{enumerate}
\item $Q$ is in regular position with respect to $P\times I$ and
\item at each regular level $P_x$, if a component $\gamma$ of $Q\cap
P_x$ is trivial in $P_x$, then $\gamma$ is also trivial in $Q$.
\end{enumerate}
\end{definition}

In this section, we  assume $Q$ is irreducible with respect to the
sweepout $P\times I$.  We first perform some isotopy on $Q$ to
eliminate center tangencies and trivial intersection curves.
\fullref{Lcircle} can be viewed as a special case of a theorem of
Thurston \cite{T} and \cite[Theorem 7.1]{G}.

\begin{lemma}\label{Lcircle}
Let $Q$ be an embedded surface in $M$ and suppose $Q$ is irreducible
with respect to the sweepout $P\times I$.  Then, one can perform an
isotopy on $Q$ so that
\begin{enumerate}
\item $Q\cap (G^U\cup G^V)$ consists of finitely many points, $Q$ is
transverse to $P^U\cup P^V$, and $Q\cap (P^U\cup P^V)$ consists of
curves essential in $Q$;
\item $Q$ is transverse to each $P_x$, $x\in (0,1)$, except for
finitely many critical levels $t_1,\dots, t_n\in (0,1)$;
\item at each critical level $t_i$, $Q$ is transverse to $P_{t_i}$
except for a saddle or circle tangency, as shown in
\fullref{Ftan}(a);
\item at each regular level $x$, every component of $Q\cap P_x$ is an
essential curve in $P_x$.
\end{enumerate}
\end{lemma}
\begin{proof}
Since $P^U\cup P^V$ is incompressible in $M$ and $M$ is irreducible,
after some standard isotopy we may assume condition (1) in the lemma
holds.

Note that the intersection of $Q$ with $P\times I$ yields a (singular)
foliation of $Q\cap M_P$ with each leaf a component of $Q\cap P_x$ for
some $x\in I$.  A singular point in the foliation is either a point in
$Q\cap (G^U\cup G^V)$ or a saddle or center tangency.

\begin{figure}
\begin{center}
\psfrag{(a)}{(a)}
\psfrag{(b)}{(b)}
\psfrag{cir}{circle (or volcano) tangency}
\psfrag{P-}{$P_{t\mp\epsilon}$}
\psfrag{P+}{$P_{t\pm\epsilon}$}
\includegraphics[width=3.5in]{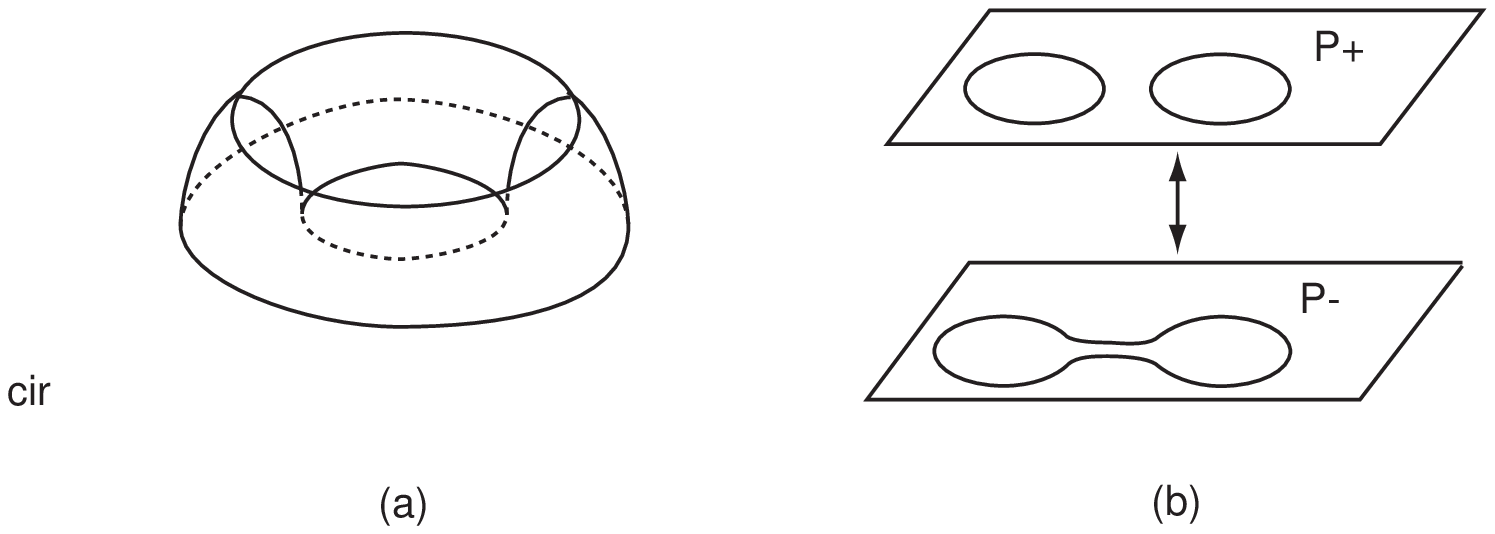}
\caption{} 
\label{Ftan}
\end{center}
\end{figure}

Let $x$ be a regular level and suppose a component $\gamma$ of $Q\cap
P_x$ is trivial in $P_x$.  Suppose $\gamma$ is innermost in $P_x$,
ie the disk bounded by $\gamma$ in $P_x$ does not contain any
other intersection curve with $Q$.  Since $Q$ is irreducible with
respect to $P\times I$, $\gamma$ bounds a disk $D_\gamma$ in $Q$.  If
the induced foliation on $D_\gamma$ contains more than one singular
point, since $\gamma$ is trivial in $P_x$, we can construct a disk
$D'\subset P\times(x-\epsilon, x+\epsilon)$ for some small $\epsilon$
such that
\begin{enumerate}
\item $\partial D'=\gamma$,
\item the induced foliation of $D'\cap (P\times I)$ consists of
parallel circles except for a singular point corresponding to a center
tangency,
\item $(Q-D_\gamma)\cup D'$ is embedded in $M$ and irreducible with
respect to $P\times I$.
\end{enumerate}
Since $M$ is irreducible, $(Q-D_\gamma)\cup D'$ is isotopic to $Q$.
Moreover, the induced foliation on $(Q-D_\gamma)\cup D'$ has fewer
singular points. So  after finitely many such operations, we may
assume that for any regular level $x$ and for any component $\gamma$
of $Q\cap P_x$ that is trivial in $P_x$, the disk bounded by $\gamma$
in $Q$ lies in $M_P$ and is transverse to $P\times (0,1)$ except for a
single center tangency.

Let $t$ be a critical level and suppose $Q\cap P_t$ contains a saddle
tangency.  Let $\epsilon$ be a sufficiently small number.  So the
component of $Q\cap (P\times [t-\epsilon, t+\epsilon])$ that contains
the saddle tangent point is a pair of pants $F$.  \fullref{Ftan}(b)
is a picture of the curves changing from $F\cap P_{t-\epsilon}$ to
$F\cap P_{t+\epsilon}$.

We first claim that at most one component of $\partial F$ is trivial
in the corresponding level surface $P_{t\pm\epsilon}$.  Let
$\gamma_1$, $\gamma_2$ and $\gamma_3$ be the 3 components of $\partial
F$ and suppose $\gamma_1$ and $\gamma_2$ are both trivial in the
corresponding level surfaces.  Then by the change of $F\cap P_x$ near
the saddle tangency as shown in \fullref{Ftan}(b), $\gamma_3$ must
also be trivial in the corresponding level surface $P_{t\pm\epsilon}$.
Since $Q$ is irreducible with respect to $P\times I$, $\gamma_1$ and
$\gamma_2$ bound disks $D_1$ and $D_2$ in $Q$ respectively.  By the
assumption above, the disk $D_i$ does not contain any saddle tangency
and hence $F\cap D_i=\gamma_i$, $i=1,2$.  Thus $F\cup D_1\cup D_2$ is
a disk in $Q$ bounded by $\gamma_3$ and $F\cup D_1\cup D_2$ contains a
saddle tangent point.  This contradicts the assumption above.  Thus at
most one component of $\partial F$ is trivial in $P_{t\pm\epsilon}$.

Let $F$ and $\gamma_i$ be as above.  Suppose $\gamma_1$ and $\gamma_2$
lie in $P_{t-\epsilon}$ and $\gamma_3$ lies in $P_{t+\epsilon}$.  If
$\gamma_1$ is trivial in $P_{t-\epsilon}$ and let $D_1$ be the disk in
$Q$ bounded by $\gamma_1$, then $F\cap D_1=\gamma_1$ as above and
$F\cup D_1$ is an annulus in $Q$ bounded by $\gamma_2\cup\gamma_3$.
Since $D_1$ is isotopic to a disk in $P_{t-\epsilon}$, we can first
push $D_1$ into $P\times [t-\epsilon, t+\epsilon]$, then as shown in
\fullref{Fiso}(a), we may perform another isotopy on $Q$ canceling
the center tangency in $D_1$ and the saddle tangency in $F$.  If
$\gamma_3$ is trivial in $P_{t+\epsilon}$, by the assumption above,
both $\gamma_1$ and $\gamma_2$ are essential in $P_{t-\epsilon}$.
Hence $\gamma_1$ and $\gamma_2$ must be parallel in $P_{t-\epsilon}$.
Let $D_3$ be the disk in $Q$ bounded by $\gamma_3$.  As above, $F\cap
D_3=\gamma_3$ and $F\cup D_3$ is an annulus in $Q$ bounded by
$\gamma_1\cup\gamma_2$.  Since $D_3$ is isotopic to the disk in
$P_{t+\epsilon}$ bounded by $\gamma_3$, we can first push the annulus
$F\cup D_3$ into a $\partial$--parallel annulus in $P\times
[t-\epsilon, t+\epsilon]$.  Then an isotopy as shown in
\fullref{Fiso}(b) can cancel the center tangency in $D_3$ and the
saddle tangency in $F$, changing $F\cup D_3$ into an annulus with a
circle (or volcano) tangency.  Note that the circle tangency is an
essential curve in the corresponding level surface $P_x$.

\begin{figure}
\begin{center}
\psfrag{(a)}{(a)}
\psfrag{(b)}{(b)}
\includegraphics[width=6cm]{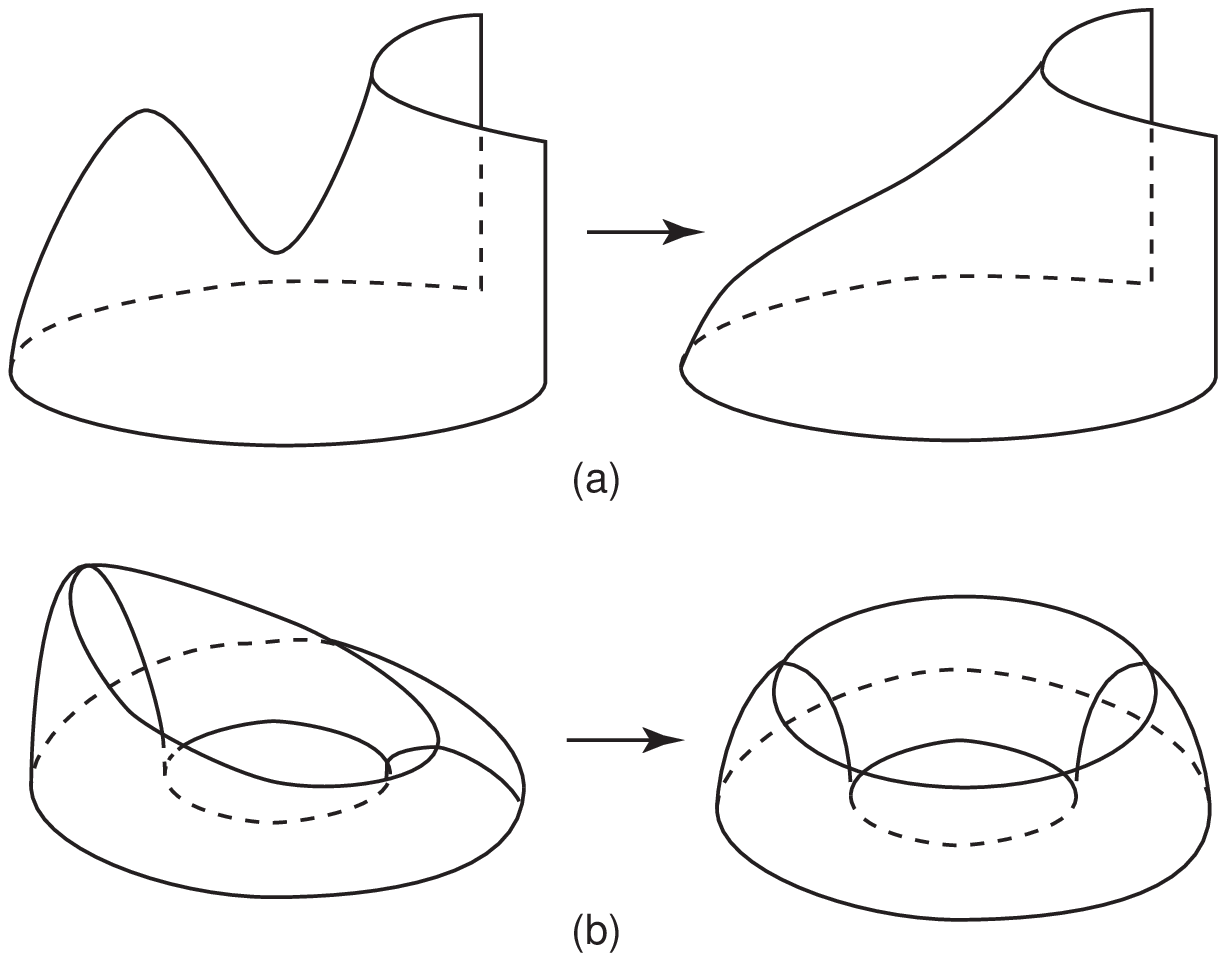}
\caption{}\label{Fiso}
\end{center}
\end{figure}

Note that condition (1) of the lemma implies that for a small
$\epsilon$, $Q\cap P_\epsilon$ and $Q\cap P_{1-\epsilon}$ consist of
essential curves in $P_\epsilon$ and $P_{1-\epsilon}$ respectively.
Since $Q$ is not a 2--sphere, a curve of $Q\cap P_x$ that is trivial in
$P_x$ will eventually meet and cancel with a saddle tangency.  Thus
after a finite number of isotopies as above, we can eliminate all the
curves of $Q\cap P_x$ that are trivial in $P_x$, and get a surface $Q$
satisfying all the conditions in the lemma.
\end{proof}

Note that a circle tangency does not create any singularity in the
foliation of $Q\cap M_P$ induced from $P\times I$.  Thus, if $Q$
satisfies the conditions in \fullref{Lcircle}, a singular point in
the foliation of $Q\cap M_P$ corresponds to either a saddle tangency
or a point in $Q\cap(G^U\cup G^V)$.  It is possible that $Q$ does not intersect $M_P=P\times I$, ie $P_t\times Q=\emptyset$ for all $t$, after isotopy.

\begin{lemma}\label{Lkey}
Let $P$ and $Q$ be as above and assume $Q$ satisfies the conditions in
\fullref{Lcircle}.  Suppose $Q\cap\Sigma_U\ne\emptyset$ and $Q\cap\Sigma_V\ne\emptyset$. Then the distance
$d(P)=d(\mathcal{U},\mathcal{V})\le 2g(Q)$.
\end{lemma}
\begin{proof} Since $Q$ is connected and $P$ is separating, $Q\cap\Sigma_U\ne\emptyset$ and $Q\cap\Sigma_V\ne\emptyset$ imply $Q\cap P_t\ne\emptyset$ for every $t$.

\begin{claim} \label{clm:claim1}
Let $t$ be a critical level and $\epsilon$
a sufficiently small number. Let $\sigma$ and $w$ be any components of
$Q\cap P_{t-\epsilon}$ and $Q\cap P_{t+\epsilon}$ respectively. Then
$d(\sigma, w)\le 1$.
\end{claim}
\begin{proof}[Proof of \fullref{clm:claim1}]
The claim is obvious if $P_t$ contains a circle tangency. So we
suppose $P_t$ contains a saddle tangency.
Let $F$ be the component of $Q\cap(P\times [t-\epsilon, t+\epsilon])$
that contains the saddle tangency. Then $F$ is a pair of pants and all
other components of $Q\cap(P\times [t-\epsilon, t+\epsilon])$ are
essential vertical annuli in $P\times[t-\epsilon, t+\epsilon]$.  If
$\sigma$ is a boundary curve of a vertical annulus, then $\sigma$ is
isotopic to a component of $Q\cap P_{t+\epsilon}$ and hence $d(\sigma,
w)\le 1$ for any curve $w$ in $Q\cap P_{t+\epsilon}$.  If neither
$\sigma$ nor $w$ is a boundary curve of a vertical annulus, then
$\sigma$ and $w$ are components of $\partial F$ and $d(\sigma, w)=1$
as shown in \fullref{Ftan}(b).
\end{proof}

Let $s_0<\cdots<s_n$ be a collection of regular levels such that
$s_0=\delta$, $s_n=1-\delta$ for a small $\delta$ and there is exactly
one saddle or circle tangency in each $P\times(s_i, s_{i+1})$.  Let
$\Gamma_i=Q\cap P_{s_i}$ for each $i$.  

Recall that $P_0=\Sigma_U=P^U\cup G^U$ and $P_1=\Sigma_V=P^V\cup G^V$.
Since $s_0=\delta$ for a small $\delta$, we may assume
$d(\mathcal{U},\Gamma_0)$ is either 0 or 1, and if
$d(\mathcal{U},\Gamma_0)=1$ then $d(\mathcal{U},\sigma)=1$ for any
component $\sigma$ of $\Gamma_0$. Similarly, $d(\mathcal{V},\Gamma_n)$
is either 0 or 1, and if $d(\mathcal{V},\Gamma_n)=1$ then
$d(\mathcal{V},w)=1$ for any component $w$ of $\Gamma_n$.

Suppose $d(\mathcal{U},\mathcal{V})>2g(Q)$ and hence
$d(\mathcal{U},\mathcal{V})>2$.  Let $k$ be the smallest integer such
that $d(\mathcal{U},\Gamma_k)\ne 0$ and $l$ the largest integer such
that $d(\Gamma_l, \mathcal{V})\ne 0$. Since
$d(\mathcal{U},\Gamma_0)\le 1$ and $d(\mathcal{V},\Gamma_n)\le 1$, by
Claim 1 above, $d(\mathcal{U},\Gamma_k)=d(\Gamma_l, \mathcal{V})=1$
and $k\le l$.  Without loss of generality, we assume $k<l$.  Next we
show that every curve in $\Gamma_k$ is essential in $Q$.  Suppose a
curve $\gamma$ in $\Gamma_k$ is trivial in $Q$ and let $D$ be the disk
bounded by $\gamma$ in $Q$.  Since $P^U$ and $P^V$ are incompressible,
we may assume $D\subset M_P$.  Since $P$ is a strongly irreducible
Heegaard surface of $M_P$, by the no-nesting lemma of Scharlemann
\cite[Lemma 2.2]{S}, $\gamma$ must bound a disk in one of the two
compression bodies, ie either $\gamma\in\mathcal{U}$ or
$\gamma\in\mathcal{V}$.  However, $\gamma\in\mathcal{U}$ contradicts
$d(\mathcal{U},\Gamma_k)\ne 0$, and $\gamma\in\mathcal{V}$ contradicts
$d(\mathcal{U},\mathcal{V})>2$.  Thus every curve in $\Gamma_k$ must
be essential in $Q$.  Similarly every curve in $\Gamma_l$ is also
essential in $Q$.

Let $Q'=Q\cap(P\times[s_k, s_l])$, and let $U'$ and $V'$ be the two
components of $M-P\times(s_k, s_l)$ containing $G^U$ and $G^V$
respectively, $F_U=Q\cap U'$ and $F_V=Q\cap V'$.  Since $\Gamma_k$ and
$\Gamma_l$ are essential in $Q$, $F_U$, $Q'$ and $F_V$ are essential
subsurfaces of $Q=F_U\cup Q'\cup F_V$.

\begin{claim} \label{clm:claim2} Let $\sigma_k$ be any component of
$\Gamma_{k}$, then $d(\sigma_k,\mathcal{U})\le 1$.
\end{claim}

\begin{proof}[Proof of \fullref{clm:claim2}]
By the definition of $k$ and the argument above, \fullref{clm:claim2} holds if
$k=0$. If $k>0$, then $d(\mathcal{U},\Gamma_{k-1})=0$ and
$d(\mathcal{U},\Gamma_{k})=1$.  Let $w$ be a component of
$\Gamma_{k-1}$ that represents a vertex in $\mathcal{U}$.  By \fullref{clm:claim1},
for any component $\sigma_k$ of $\Gamma_{k}$,
$d(\sigma_k,\mathcal{U})\le d(\sigma_k,w)\le 1$.
\end{proof}

\begin{claim} \label{clm:claim3}  There is a component $\sigma_k$ of
$\Gamma_k$ and a component $\sigma_l$ of $\Gamma_l$ such that
$d(\sigma_k,\sigma_l)\le-\chi(Q')$.
\end{claim}

\begin{proof}[Proof of \fullref{clm:claim3}]
Let $t_1<\cdots<t_N$ be the levels in $(s_k,s_l)$ that contain the
saddle tangencies.  For a sufficiently small $\epsilon$,
$P\times[t_i+\epsilon, t_{i+1}-\epsilon]$ contains no saddle tangency
for each $i$ (to simplify notation we set $t_0+\epsilon=s_k$ and
$t_{N+1}-\epsilon=s_l$). So by the conditions in \fullref{Lcircle},
$Q\cap (P\times[t_i+\epsilon, t_{i+1}-\epsilon])$ consists of annuli
for each $i=0,\dots, N$. If $Q\cap (P\times[t_i+\epsilon,
t_{i+1}-\epsilon])$ consists of $\partial$--parallel annuli, then $Q\cap P_t=\emptyset$ for some $t$ after isotopy, a contradiction to our assumption at the beginning. Thus an annulus component $A_i$ of
$Q\cap (P\times[t_i+\epsilon, t_{i+1}-\epsilon])$ is vertical.  We
choose $\gamma_i$ to be a meridian circle in $A_i$ for each $i$ and
assume $\sigma_k=\gamma_0=A_0\cap P_{s_k}\subset\Gamma_k$ and
$\sigma_l=\gamma_N=A_N\cap P_{s_l}\subset\Gamma_l$.  Since each $A_i$
is vertical, $\gamma_i$ is parallel to a component of $Q\cap
P_{t_{i+1}-\epsilon}$. Similarly $\gamma_{i+1}$ is parallel to a
component of  $Q\cap P_{t_{i+1}+\epsilon}$.  By \fullref{clm:claim1}, $d(\gamma_i,
\gamma_{i+1})\le 1$ and hence $d(\sigma_k,\sigma_l)=d(\gamma_0,
\gamma_{N})\le N$.  Moreover, since the only singular points in the
induced foliation of $Q'$ are the saddle tangencies, by a standard
index argument, $-\chi(Q')=N$ and hence
$d(\sigma_k,\sigma_l)\le-\chi(Q')$.
\end{proof}

Since $Q'$, $F_U$ and $F_V$ are essential subsurfaces of $Q$,
$\chi(Q')\ge\chi(Q)$.  By \fullref{clm:claim2}, $d(\sigma_k,\mathcal{U})\le 1$ and
similarly $d(\sigma_l, \mathcal{V})\le 1$.  Therefore,
$d(\mathcal{U},\mathcal{V})\le
d(\mathcal{U},\sigma_k)+d(\sigma_k,\sigma_l)+d(\sigma_l,
\mathcal{V})\le 1-\chi(Q')+1\le 2-\chi(Q)=2g(Q)$.
\end{proof}

\begin{corollary}\label{Cmain}
Let $P$ and $Q$ be as in \fullref{Tmain}.  Then
\fullref{Tmain} holds if $Q$ is incompressible.
\end{corollary}
\begin{proof}
If $Q$ is incompressible, then $Q$ can be isotoped to be irreducible
with respect to $P\times I$. Moreover, if $Q\cap\Sigma_U=\emptyset$, then since $Q$ is incompressible, $Q$ can be isotoped out of the compression body $M_P-N(\Sigma_U)$. Hence $Q\cap M_P=\emptyset$ after isotopy and part (2) of \fullref{Tmain} holds.  Now \fullref{Cmain} follows from
 \fullref{Lcircle} and \fullref{Lkey}.
\end{proof}

\section{The graphics of sweepouts}

In this section, we suppose $Q$ is separating, bicompressible and
strongly irreducible.

Let $X$ and $Y$ be the closure of the 2 components of $M-Q$.  Let
$Q^X$ and $Q^Y$ be the possibly disconnected surfaces obtained by
maximal compressing $Q$ in $X$ and $Y$ respectively and capping off
2--sphere components by 3--balls.  Similar to the argument on $P^U$ and
$P^V$ above, we may assume $Q^X\subset \Int(X)$ and $Q^Y\subset \Int(Y)$
are incompressible in $M$.    Furthermore, $Q^X\cup Q^Y$ bounds a
submanifold $M_Q$ of $M$ and $Q$ is a strongly irreducible Heegaard
surface of $M_Q$.  If $X$ is a handlebody, then $Q^X=\emptyset$.  If
$Q$ is a Heegaard surface of $M$, we may view $M_Q=M$.

As in \fullref{Skey}, the surface $Q$ cuts $M_Q$ into a pair of
compression bodies $X\cap M_Q$ and $Y\cap M_Q$.  Let graphs
$G^X\subset X\cap M_Q$ and $G^Y\subset Y\cap M_Q$ be the spines of the
two compression bodies and let $\Sigma_X=Q^X\cup G^X$ and
$\Sigma_Y=Q^Y\cup G^Y$. Then $M_Q-(\Sigma_X\cup\Sigma_Y)$ is
homeomorphic to $Q\times (0,1)$.

Now we consider the two sweepouts $H \co P\times (I,\partial I)\to
(M_P,\Sigma_U\cup\Sigma_V)$ and $H' \co Q\times (I,\partial I)\to
(M_Q,\Sigma_X\cup\Sigma_Y)$.  Let $P_t=H(P\times\{t\})$ and
$Q_t=H'(Q\times\{t\})$, $t\in I$. We may assume $Q_0=\Sigma_X$,
$Q_1=\Sigma_Y$ and $Q_x$ is isotopic to $Q$ for each $x\in(0,1)$.

The graphic $\Lambda$ of the sweepouts, defined in \cite{RS}, is the
set of points $(s,t)\in(0,1)\times (0,1)$ such that $P_s$ is not
transverse to $Q_t$.  We briefly describe the graphic below and refer
to \cite{RS} for more details.  As in \cite{RS}, Cerf theory implies
that after some isotopy, we may assume that $\Lambda$ is a graph in
$(0,1)\times(0,1)$ whose edges are the set of points $(s,t)$ for which
$P_s$ is transverse to $Q_t$ except for a single saddle or center
tangency.  There are two types of vertices in $\Lambda$,
birth-and-death vertices and crossing vertices, as shown in
\fullref{Fver}(a).  Moreover, each arc $(0,1)\times\{x\}$ contains
at most one vertex, $x\in(0,1)$. The complement of $\Lambda$,
$(0,1)\times(0,1)-\Lambda$, is a finite collection of regions.  Note
that for every $(s,t)$ in $(0,1)\times(0,1)-\Lambda$, $P_s$ is
transverse to $Q_t$, and for any two points $(s,t)$ and $(s',t')$ in
the same region, $P_s\cap Q_t$ and $P_{s'}\cap Q_{t'}$ have the same
intersection pattern.

\begin{figure}
\begin{center}
\psfrag{(a)}{(a)}
\psfrag{(b)}{(b)}
\psfrag{(c)}{(c)}
\psfrag{(d)}{(d)}
\psfrag{X}{$X$}
\psfrag{Y}{$Y$}
\psfrag{Z}{$Z$}
\psfrag{b}{$b$}
\psfrag{b'}{$b\pm\epsilon$}
\includegraphics[width=4.5in]{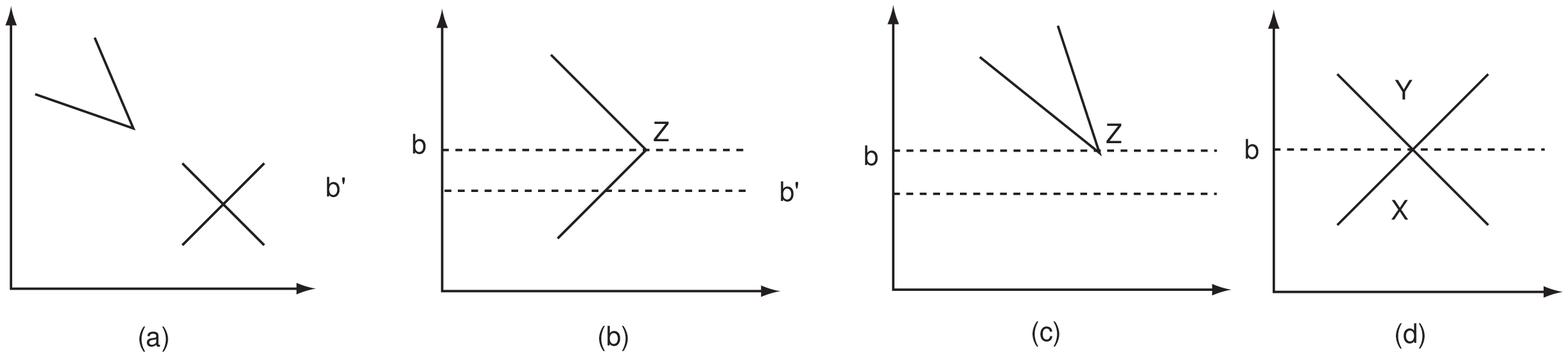}
\caption{}\label{Fver}
\end{center}
\end{figure}

Let $(s,t)\in (0,1)\times(0,1)-\Lambda$.
Suppose there are disks or annuli $D_P\subset P_s$ and $D_Q\subset Q_t$ with $D_P\cap D_Q=\partial D_P=\partial D_Q\subset P_s\cap Q_t$.  Suppose $D_P$ is parallel to $D_Q$ (fixing $\partial D_P=\partial D_Q$) in $M$ and suppose $D_P\cup D_Q$ bounds a 3-ball or solid torus $E$.  Moreover, suppose $Q_t\cap E=D_Q$. 
Then we can perform an isotopy on $Q_t$ by pushing 
$D_Q$ across $E$ and remove the intersection $\partial D_P=\partial D_Q$. This isotopy is the same as the operation that changes $Q_t$ to $(Q_t-D_Q)\cup D_P$ and then perturbs the resulting surface.  We call such an isotopy a \emph{trivial isotopy} on $Q_t$ at $P_s$. 
We may view a trivial isotopy on $Q_t$ is associated with the disk or annulus $D_Q\subset Q_t$.  Suppose we are to perform another trivial isotopy on $Q_t$ at $P_{s'}$ and let $D_Q'\subset Q_t$ be the disk or annulus in the isotopy as above.  Then $D_Q$ and $D_Q'$ are either disjoint or nested in $Q_t$. Thus either the two trivial isotopies are disjoint or we can view one isotopy as a middle step of the other.

\begin{labelling} For any $Q_t$, we use $X_t$ (resp. $Y_t$) to denote
the component of $M-Q_t$ that contains $\Sigma_X$ 
(resp. $\Sigma_Y$).  
We label a region, ie a component of
$(0,1)\times(0,1)-\Lambda$, $X$ (resp. $Y$) if for a point $(s,t)$ in
the region, either (1) there is a component of $P_s\cap Q_t$ that is trivial in
$P_s$ but bounds an essential disk in $X_t$ (resp. $Y_t$), or (2) $\Sigma_U$ or
$\Sigma_V$ lies in $Y_t$ (resp. $X_t$) after some \emph{trivial isotopies} on $Q_t$ at finitely many regular levels $P_x$.  We label $t\in (0,1)$ $X$ (resp. $Y$) if the horizontal line segment 
$(0,1)\times\{t\}$ intersects a region labelled $X$ (resp. $Y$).  Note that since a trivial isotopy does not increase $|\Sigma_U \cap Q_t|$ or$|\Sigma_V \cap Q_t|$,
if $t$ is not labelled, $Q_t\cap\Sigma_U\ne\emptyset$ and $Q_t\cap\Sigma_V\ne\emptyset$ after any trivial isotopies.
\end{labelling}

\begin{lemma}\label{Lsmall}
Either \fullref{Tmain} holds or for a sufficiently small
$\delta>0$, $\delta$ is labelled $X$ and $1-\delta$ is labelled $Y$.
\end{lemma}
\begin{proof}
For a sufficiently small $\delta>0$, $H'(Q\times [0,\delta])$ is a
small neighborhood of $\Sigma_X=Q^X\cup G^X$.  If $P_s\cap
G^X\ne\emptyset$ for some $s$, then by definition, $\delta$ is labelled
$X$ for a sufficiently small $\delta$.  Suppose $\delta$ is not
labelled $X$, then the graph $G^X$ must be disjoint from $M_P=H(P\times
I)$.  Moreover, if $Q^X\cap P_t=\emptyset$ for some $t$ after isotopy,
since $Q^X$ is incompressible, we can isotope $Q^X$ out of the two
compression bodies $M_P-P_t$. Hence, $Q_\delta\cap M_P=\emptyset$ 
after isotopy and part (2) of \fullref{Tmain} holds.  
If $Q^X\cap P_t\ne\emptyset$ for all $t$, since $Q^X$ is incompressible, 
by \fullref{Cmain}, $d(P)\le 2g(Q^X)\le 2g(Q)$ and \fullref{Tmain} 
follows.  The proof for $1-\delta$ is similar.
\end{proof}

\begin{lemma}\label{Lboth}
Either \fullref{Tmain} holds or no $t\in (0,1)$ is labelled both $X$ and $Y$.
\end{lemma}
\begin{proof}
We first remark that if $\Sigma_U\subset Y_t$ then one cannot move $\Sigma_U$ to $X_t$ by a trivial isotopy, since if this happens, then $\Sigma_U$ must lie in $E$, where $E$ is the 3-ball or solid torus in the trivial isotopy described above.  However, since $g(P)\ge 2$ and $P$ is strongly irreducible, $\Sigma_U$ cannot lie in a 3-ball or solid torus by \cite{CG}.  So by our labelling, if $t$ is labelled both $X$ and $Y$, then one can always find $s_1\ne s_2$ such that $(s_1,t)$ and $(s_2,t)$ are points in $(0,1)\times(0,1)-\Lambda$ and one of the following 3 cases occurs.

(1).  A component of $P_{s_1}\cap Q_t$ contains a
curve bounding an essential disk $D_X$ in 
$X_t$ and a component of $P_{s_2}\cap Q_t$ contains a
curve bounding an essential disk $D_Y$ in 
$Y_t$. In this case, since $s_1\ne s_2$, $\partial D_X\cap\partial
D_Y=\emptyset$ in $Q_t$, which contradicts the assumption that $Q$ is
strongly irreducible.

(2). After trivial isotopies, $\Sigma_U\subset Y_t$ and $\Sigma_V\subset X_t$.  This means that $Q_t\subset P\times(0,1)\subset M_P$ and $Q_t$ separates $\Sigma_U$ and $\Sigma_V$ in $M_P$.  If $Q_t$ is incompressible in $P\times (0,1)$, then $Q_t$ is isotopic to $P$ and \fullref{Tmain} holds.  Suppose $Q_t$ is compressible in $P\times (0,1)$.  Similar to the construction of $M_Q$ earlier, by maximal compressing $Q_t$ in $P\times (0,1)$ on both sides and capping off 2-sphere components, we obtain a submanifold $M_Q'$ of $P\times(0,1)$ such that $Q_t$ is a strongly irreducible Heegaard surface of $M_Q'$.  Moreover, by \cite{CG}, $\partial M_Q'$ is incompressible in $P\times(0,1)$.  So each component of $\partial M_Q'$ is parallel to $P$ and $M_Q'$ must be a product of $P$ and an interval.  Thus we can view $Q_t$ as a strongly irreducible Heegaard surface of a product $P\times [0,1]$. By Scharlemann--Thompson \cite{ST1},  either $Q_t$ is isotopic to $P$ or $Q_t$ cuts $P\times [0,1]$ into a handlebody and a compression body. In the later case, both $\Sigma_U$ and $\Sigma_V$ lie in $Y_t$ (or both in $X_t$), a contradiction.

(3). After trivial isotopies, $\Sigma_U\subset Y_t$ and a component of $P_{s_1}\cap Q_t$ contains a curve $\gamma$ that is trivial in $P_{s_1}$ and bounds an essential disk $D$ in $Y_t$.  Note that if a component of $P_{s_1}\cap Q_t$ also bounds an essential disk in $X_t$, then this contradicts that $Q$ is strongly irreducible as in case (1). Thus, after some isotopy on $Q_t$, we may assume that $\gamma$ is innermost in $P_{s_1}$ and the disk $D$ bounded by $\gamma$ in $P_{s_1}$ is an essential disk in $Y_t$.   Since $\Sigma_U\subset Y_t$ and $D\subset Y_t-\Sigma_U$, by maximal compressing $Q_t$ in $Y_t-\Sigma_U$  and capping off 2-sphere components, we obtained a (possibly disconnected) surface $Q_t^Y$. Note that $Q_t^Y\ne\emptyset$ because $\Sigma_U$ is not contained in a 3-ball.  Since $Q_t$ is strongly irreducible, by \cite{CG}, $Q_t^Y$ is incompressible in $M-\Sigma_U$. 
Note that if $P$ is a Heegaard surface of a closed manifold $M$, this is already a contradiction since $Q_t^Y$ lies in the handlebody $M-N(\Sigma_U)$ and cannot be incompressible. $Q_t^Y$ cuts $Y_t$ into $H_1$ and $H_2$, where $H_2$ is the compression body bounded by $Q_t$ and $Q_t^Y$.  Since the compressions on $Q_t$ are disjoint from $\Sigma_U$ and since $\Sigma_U$ does not lie in a 3-ball, $\Sigma_U\cap H_2=\emptyset$.  Hence $\Sigma_U\subset H_1$.  Since $Q_t^Y$ is incompressible in $M-\Sigma_U$, we can push $Q_t^Y$ out of the compression body $M_P-N(\Sigma_U)$ or equivalently push $M_P-N(\Sigma_U)$ into $H_1$.  So we can isotope $M_P$ into $H_1$. In particular, $Q_t\cap M_P=\emptyset$ after isotopy and part(2) of \fullref{Tmain} holds. 
\end{proof}

\begin{lemma}\label{Lnon}
If $t\in (0,1)$ has no label and $(0,1)\times\{t\}$ contains no vertex
of $\Lambda$, then $Q_t$ is irreducible with respect to $P\times I$
and \fullref{Tmain} holds.
\end{lemma}
\begin{proof}
Since $(0,1)\times\{t\}$ contains no vertex of $\Lambda$, $Q_t$ is in
regular position with respect to $P\times I$.   For any
$(s,t)\notin\Lambda$, suppose a curve $\gamma$ in $P_s\cap Q_t$ is
trivial in $P_s$.  If $\gamma$ is an essential curve in $Q_t$, by assuming $\gamma$ to be an innermost such curve, the disk bounded by $\gamma$ in $P_s$ can be isotoped to be an essential disk in either $X_t$ or $Y_t$.
Since $t\in
(0,1)$ has no label, $\gamma$ must be trivial in $Q_t$.  Thus by
definition, $Q_t$ is irreducible with respect to $P\times I$. So after
isotopy we may assume $Q$ satisfies the conditions in \fullref{Lcircle}. 
Moreover, since $t$ has no label, $Q_t\cap\Sigma_U\ne\emptyset$ and 
$Q_t\cap\Sigma_V\ne\emptyset$ after the isotopy in the proof of 
\fullref{Lcircle}. So \fullref{Tmain} follows from \fullref{Lkey}.
\end{proof}

Suppose \fullref{Tmain} is not true. Then by \fullref{Lsmall},
for a small $\delta$, $\delta$ is labelled $X$ and $1-\delta$ is
labelled $Y$.  As $t$ changes from $\delta$ to $1-\delta$, the label
changes from $X$ to $Y$.  So by \fullref{Lboth} and \fullref{Lnon},
there must be a number $b\in (0,1)$ such that
\begin{enumerate}
\item $(0,1)\times\{b\}$ contains a vertex of $\Lambda$ and
\item $b$ has no label and
\item $b-\epsilon$ is labelled $X$ and $b+\epsilon$ is labelled $Y$ for
sufficiently small $\epsilon>0$.
\end{enumerate}
Let $Z=(a,b)$ be the vertex of $\Lambda$ in $(0,1)\times\{b\}$.  If
$Z$ is a birth-and-death vertex, then since no region that intersects
$(0,1)\times\{b\}$ is labelled, as shown in \fullref{Fver}(b) and
(c), after perturb $(0,1)\times\{b\}$ a little, we can find a line
segment $(0,1)\times\{b\pm\epsilon\}$ that does not intersect any
labelled region, a contradiction to our assumption above.  Therefore,
$Z=(a,b)$ must be a crossing vertex.  \fullref{Fver}(d) is a picture of $Z$.

Since $Z=(a,b)$ is a crossing vertex, as explained in \cite{RS} (see
Kobayashi--Saeki \cite[Figure 2.6]{KS}), $P_a$ is transverse to $Q_b$ except for two
saddle tangencies.   Since $b$ is not labelled, for any $s\ne a$ in
$(0,1)$, either (1) $P_s\cap Q_b$ contains a single center or saddle
tangency or (2) $P_s$ is transverse to $Q_b$ and if a component of
$P_s\cap Q_b$ is trivial in $P_s$ then it is also trivial in $Q_b$.  Moreover, after  trivial isotopies, $Q_b\cap\Sigma_U\ne\emptyset$ and $Q_b\cap\Sigma_V\ne\emptyset$. Since $P$ is separating and $Q$ is connected, this implies that $Q_b\cap P_s\ne\emptyset$ for all $s\in I$.

Now we consider $Q_b\cap (P\times[a-\epsilon, a+\epsilon])$ for a
small $\epsilon$.  Let $F$ be the union of the components of $Q_b\cap
(P\times[a-\epsilon, a+\epsilon])$ that contain the two saddle
tangencies.  Thus $F$ is either the union of two disjoint pairs of
pants or a connected surface with $\chi(F)=-2$.  All other components
of $Q_b\cap (P\times[a-\epsilon, a+\epsilon])$, denoted by $A_1,\dots,
A_m$, are vertical annuli in $P\times[a-\epsilon, a+\epsilon]$.

Next we consider the case that a component of $Q_b\cap
P_{a\pm\epsilon}$ is trivial in $P_{a\pm\epsilon}$.  If a component
$\gamma$ of $\partial A_i$, $i=1,\dots, m$, is trivial and innermost
in $P_{a\pm\epsilon}$, then by our assumption, $\gamma$ bounds a disk
$D_\gamma$ in $Q_b$.  We can perform a trivial isotopy on $Q_b$ by
pushing the disk $D_\gamma\cup A_i$ away from $P\times[a-\epsilon,
a+\epsilon]$.  Thus, after a finite number of such operations, we may
assume the boundary of every annular component $A_i$ is essential in
$P_{a\pm\epsilon}$.

Suppose a component $\gamma$ of $\partial F$ is an innermost trivial
curve in $P_{a\pm\epsilon}$.  So $\gamma$ bounds a disk $D_\gamma$ in
$Q_b$.  If $D_\gamma$ contains a component of $F$, then as in the
proof of \fullref{Lcircle}, after replacing $D_\gamma$ by a disk
which is transverse to every $P_x$ except for a single center
tangency, we get a surface isotopic to $Q_b$ and has at most one
saddle tangency in $P\times[a-\epsilon, a+\epsilon]$.  This means that
after the isotopy, $Q_b$ is irreducible with respect to $P\times I$
and \fullref{Tmain} follows from \fullref{Lkey} and \fullref{Lnon}.  So we may
assume that $D_\gamma\cap F=\gamma$ for any component $\gamma$ of
$\partial F$ that is trivial in $P_{a\pm\epsilon}$.

Let $\hat{F}$ be the union of $F$ and all the disks $D_\gamma$ in
$Q_b$ bounded by $\partial F$ as above.   We may push all such disks
$D_\gamma$ into $P\times(a-\epsilon, a+\epsilon)$ and isotope
$\hat{F}$ into a surface properly embedded in $P\times[a-\epsilon,
a+\epsilon]$.  By the construction, $\partial\hat{F}$ is essential in
$P_{a\pm\epsilon}$.  So $\hat{F}$ has no disk component.  If $\hat{F}$
consists of annuli, then since $\partial\hat{F}$ is essential in
$P_{a\pm\epsilon}$, each annulus is either vertical or
$\partial$--parallel in $P\times[a-\epsilon, a+\epsilon]$.  Thus, after
some isotopy, $Q_b$ becomes irreducible with respect to $P\times I$
and \fullref{Tmain} follows from \fullref{Lkey} and \fullref{Lnon}. So we may assume $\chi(\hat{F})$ is
either $-2$ or $-1$, ie at most one component of $\partial F$ is
trivial in $P_{a\pm\epsilon}$.

Suppose $\chi(\hat{F})=-1$. If $\hat{F}$ is a once-punctured torus,
then $\hat{F}$ must be incompressible in $P\times[a-\epsilon,
a+\epsilon]$. Otherwise a compression on $\hat{F}$ yields a disk,
contradicting that $\partial\hat{F}$ is essential in
$P_{a\pm\epsilon}$.  As $\hat{F}$ is properly embedded in the product
$P\times[a-\epsilon, a+\epsilon]$, $\hat{F}$ must be
$\partial$--compressible.  A $\partial$--compression on $\hat{F}$ yields
an incompressible annulus with both boundary circles in
$P_{a-\epsilon}$ (or $P_{a+\epsilon}$). So the resulting annulus is
$\partial$--parallel.  Since $\hat{F}$ is incompressible, this implies
that $\hat{F}$ itself is $\partial$--parallel.  Hence we can perform an
isotopy on $\hat{F}$ so that $Q_b$ becomes irreducible with respect to
$P\times I$.  Similarly, if $\hat{F}$ is a pair of pants, then
$\hat{F}$ must be incompressible but $\partial$--compressible. So a
$\partial$--compression on $\hat{F}$ yields one or two incompressible
annuli.  This implies that either $\hat{F}$ is $\partial$--parallel or
we can perform an isotopy on $\hat{F}$ so that $\hat{F}$ is transverse
to each $P_x$ except for a single saddle tangency. In either case, we
can isotope $\hat{F}$ so that $Q_b$ becomes irreducible with respect
to $P\times I$ and \fullref{Tmain} follows from \fullref{Lkey} and \fullref{Lnon}.

Therefore, we may assume $\chi(\hat{F})=-2$. Hence $F=\hat{F}$ and
every component of $\partial F$ is essential in $P_{a\pm\epsilon}$.

Since $b$ is not labelled and since every component of $\partial F$
above is essential in $P_{a\pm\epsilon}$, at each regular level
$x\in(0,1)$, if a component of $P_x\cap Q_b$ is trivial in $P_x$, then
it must also be trivial in $Q_b$.   Thus, we can apply
\fullref{Lcircle} on $Q_b\cap(P\times ([0, a-\epsilon]\cup
[a+\epsilon,1]))$. So after some isotopies, $Q_b$ satisfies all the
conditions in \fullref{Lcircle} except for the level $P_a$ where
$P_a\cap Q_b$ contains 2 saddle tangencies.  Moreover, since $b$ is not labelled, $Q_b\cap\Sigma_U\ne\emptyset$ and $Q_b\cap\Sigma_V\ne\emptyset$. Hence $Q_b\cap P_s\ne\emptyset$ for every $s$.

\begin{claima} \label{clm:claimA} Let $\sigma$ and $w$ be any components of
$Q_b\cap P_{a-\epsilon}$ and $Q_b\cap P_{a+\epsilon}$ respectively.
Then $d(\sigma,w)\le 2=-\chi(F)=-\chi(Q_b\cap(P\times[a-\epsilon,
a+\epsilon]))$.
\end{claima}

\begin{proof}[Proof of Claim A]
If $\sigma$ is a boundary curve of a vertical annulus component of
$Q_b\cap(P\times[a-\epsilon, a+\epsilon])$, then $\sigma$ is isotopic
to a component of $Q\cap P_{a+\epsilon}$ and hence $d(\sigma, w)\le 1$
for any curve $w$ in $Q\cap P_{a+\epsilon}$.  So we may assume neither
$\sigma$ nor $w$ is a boundary curve of a vertical annulus.  Thus
$\sigma$ and $w$ are both components of $\partial F$.

Let $\Omega$ be the union of the components of $P_a\cap Q_b$ that
contain the 2 saddle tangent points.  So $\Omega$ is a possibly
disconnected graph with 2 vertices of valence 4.  Let $N(\Omega)$ be a
regular neighborhood of $\Omega$ in $P_a$ and let $\pi \co P\times I\to
P_a$ be the projection, then $\pi(\partial F)\subset N(\Omega)$ after
isotopy.  Since $P$ has genus at least 2, there must be an essential
curve $\alpha$ in $P_a$ disjoint from  $N(\Omega)$.  So  $d(\sigma,
w)\le d(\sigma,\alpha)+ d(\alpha, w)\le 2=-\chi(F)$.
\end{proof}

Now \fullref{Tmain} follows from the argument in the proof of
\fullref{Lkey}.  As in the proof of \fullref{Lkey}, let
$s_0<\cdots <s_n$ be a collection of regular levels such that
$s_0=\delta$, $s_n=1-\delta$ for a small $\delta$ and there is exactly
one critical level in each $P\times(s_i, s_{i+1})$.  Let
$\Gamma_i=Q\cap P_{s_i}$ for each $i$.

Since we assume $Q$ is bicompressible in this section and since $M$ is
irreducible, if $Q$ is a torus, then $M$ must be a lens space and $P$
and $Q$ must be isotopic Heegaard surfaces of the lens space
(see Bonahon--Otal \cite{BO}). So we may assume $g(Q)\ge 2$.

Suppose $d(\mathcal{U},\mathcal{V})>g(Q)$. Since $g(Q)\ge 2$, we have
$d(\mathcal{U},\mathcal{V})>4$. Let $k$ be the smallest integer such
that $d(\mathcal{U},\Gamma_k)\ne 0$ and $l$  the largest integer such
that $d(\Gamma_l, \mathcal{V})\ne 0$.  By Claim A above and 
\fullref{clm:claim1} in
the proof of \fullref{Lkey}, $d(\mathcal{U},\Gamma_k)$ and
$d(\Gamma_l, \mathcal{V})$ are either 1 or 2 and $k\le l$.  Without
loss of generality, we assume $k<l$.  

Similar to the proof of \fullref{Lkey}, $\Gamma_k$ and $\Gamma_l$
must be essential in $Q_b$.  Let $Q'=Q_b\cap(P\times[s_k, s_l])$, and
let $U'$ and $V'$ be the two components of $M-P\times(s_k, s_l)$
containing $G^U$ and $G^V$ respectively, $F_U=Q_b\cap U'$ and
$F_V=Q_b\cap V'$.  Since $\Gamma_k$ and $\Gamma_l$ are essential in
$Q_b$, $F_U$, $Q'$ and $F_V$ are essential subsurfaces of $Q_b=F_U\cup
Q'\cup F_V$.

\begin{claimb} \label{clm:claimB} Let $\sigma_k$ be any component of
$\Gamma_{k}$, then $d(\sigma_k,\mathcal{U})\le 1-\chi(F_U)$.
\end{claimb}

\begin{proof}[Proof of Claim B]
If a component $A$ of $F_U$ is a $\partial$--parallel annulus in $U'$,
then we may first isotope $A$ into  $P\times(s_k-\epsilon,s_k]$.  Then
we  isotope $A$ so that $A$ is transverse to each $P_x$ except for a
circle  tangency.  Since $\partial F_U$ is essential in $P_{s_k}$,
after the isotopy, $Q_b$ still satisfies the conditions in
\fullref{Lcircle} except at the level $P_a$ as above.  Now we push
$A$ out of $U'$.  After the
isotopy, we still have $d(\mathcal{U},\Gamma_k)\ne 0$.  If $k$ is no
longer the smallest number so that $d(\mathcal{U},\Gamma_k)\ne 0$
after the isotopy, then we can find a new $k$ and proceed as above.
Eventually $F_U$ does not contain any $\partial$--parallel annulus
after some isotopies.  We can view these isotopies as trivial isotopies, so by our assumptions above, $Q_b\cap\Sigma_U\ne\emptyset$ after the isotopies.

We first show that $d(\sigma_k,\mathcal{U})\le 2$.  As in the proof of
\fullref{Lkey}, $d(\sigma_k,\mathcal{U})\le 1$ if $k=0$.  So we may
assume $k>0$. By the definition of $k$,
$d(\mathcal{U},\Gamma_{k-1})=0$.  Thus there is a component $w$ of
$\Gamma_{k-1}$ representing a vertex in $\mathcal{U}$.  By Claim A
above and the \fullref{clm:claim1} in the proof of \fullref{Lkey}, $d(\sigma_k,
w)\le 2$ and hence $d(\sigma_k,\mathcal{U})\le 2$.

Since $F_U$ is an essential subsurface of $Q_b$, $\chi(F_U)\le 0$.
Since $d(\sigma_k,\mathcal{U})\le 2$ and $\chi(F_U)\le 0$, to prove
the claim, we only need to consider the case that $\chi(F_U)=0$.
Suppose $\chi(F_U)=0$. Since $d(\mathcal{U},\Gamma_k)\ne 0$, $F_U$
consists of incompressible annuli in $U'$.  Let $A$ be the component
of $F_U$ that contains $\sigma_k$.  If $A$ is also
$\partial$--incompressible, then $A$ can be isotoped away from any
compressing disk of $U'$ and hence $d(\sigma_k,\mathcal{U})\le
1=1-\chi(F_U)$. If $A$ is $\partial$--compressible, then since $F_U$
contains no $\partial$--parallel annulus, a $\partial$--compression on
$A$ yields a compressing disk of $U'$ disjoint from $A$. Thus,
$d(\sigma_k,\mathcal{U})\le 1=1-\chi(F_U)$ in any case.
\end{proof}

Similar to Claim B, for any component $\sigma_l$ of $\Gamma_l$,
$d(\mathcal{V}, \sigma_l)\le 1-\chi(F_V)$. Although $P_a\cap Q_b$
contains 2 saddle tangencies, by Claim A and our assumptions on $Q_b$,
\fullref{clm:claim3} in the proof of \fullref{Lkey} also holds in this case,
ie there is a component $\sigma_k$ of $\Gamma_k$ and a component
$\sigma_l$ of $\Gamma_l$ such that $d(\sigma_k,\sigma_l)\le-\chi(Q')$.

Since $Q'$, $F_U$ and $F_V$ are essential subsurfaces of $Q_b$,
$d(\mathcal{U},\mathcal{V})\le
d(\mathcal{U},\sigma_k)+d(\sigma_k,\sigma_l)+d(\sigma_l,
\mathcal{V})\le 1-\chi(F_U)-\chi(Q')+1-\chi(F_V)=2-\chi(Q)=2g(Q)$.
Thus \fullref{Tmain} is proved. \hfill \qed

\end{psfrags}
\end{document}